\documentclass[10pt]{amsart}
\usepackage{enumerate,amsmath,amssymb,latexsym,
amsfonts, amsthm, amscd, MnSymbol}

\theoremstyle{plain}

\numberwithin{equation}{section}

\addtocounter{section}{-1}

\begin{document}

\title {Neville's primitive elliptic functions: the case ${\bf g_3 = 0}$}

\date{}

\author[P.L. Robinson]{P.L. Robinson}

\address{Department of Mathematics \\ University of Florida \\ Gainesville FL 32611  USA }

\email[]{paulr@ufl.edu}

 \keywords{}

\begin{abstract}

The vanishing of the invariant $g_3$ attached to a lattice $\Lambda$ singles out a midpoint lattice and yields a square-root of the associated Weierstrass function $\wp_{\Lambda}$. 

\end{abstract}

\maketitle

\medbreak

Neville's {\it Jacobian Elliptic Functions} [1] is a peerless classic in its field. It is therefore with some reticence that we draw attention to a minor oversight in its presentation of the primitive functions as meromorphic square-roots of the shifted Weierstrass function $\wp$. 

\medbreak 

The oversight occurs on page 50 of [1]: there it is stated that `The zeros of $\wp z$ are simple, and the branches of $(\wp z)^{\frac{1}{2}}$ can not be separated'. This is not quite correct: if we denote by $\Lambda$ the period/pole lattice of $\wp$ then the zeros of $\wp$ are simple {\it except} in case the invariant $g_3 (\Lambda)$ is zero. We note that the statement asserting simplicity of the zeros of $\wp$ is not made in the prequel [2]. 

\medbreak 

In a little more detail, let $\Lambda \subset \mathbb{C}$ be any lattice; the associated Weierstrass function $\wp = \wp_{\Lambda}$ is then defined by 
$$\wp (z) = z^{-2} + \sum_{0 \neq \lambda \in \Lambda} \{ (z - \lambda)^{-2} - \lambda^{-2} \}$$
and has $\Lambda$ as both its period lattice and its pole lattice. The zeros of the derivative $\wp '$ are precisely those $z \notin \Lambda$ such that $2 z \in \Lambda$ and they make up three congruent lattices. If 
$$\Lambda = \{ 2n_1 \omega_1 + 2n_2 \omega_2 : n_1, n_2 \in \mathbb{Z} \}$$
 and $\omega_1 + \omega_2 + \omega_3 = 0$ then these three {\it midpoint lattices} are $\omega_1 + \Lambda$, $\omega_2 + \Lambda$, $\omega_3 + \Lambda$; the values of $\wp$ at each point of these lattices are denoted by $e_1, e_2, e_3$ respectively. Among their many properties, these distinct {\it midpoint constants} satisfy 
$$e_1 + e_2 + e_3 = 0$$
and 
$$e_1 e_2 e_3 = g_3/4$$
where the invariant $g_3 = g_3 (\Lambda)$ is defined by
$$g_3 = 140 \sum_{0 \neq \lambda \in \Lambda} \lambda^{-6}.$$

\medbreak 

It follows at once that if $g_3 = 0$ then precisely one of the midpoint constants vanishes: say $0 = e_p = \wp (\omega_p)$; as $\wp ' (\omega_p) = 0$ also, $\omega_p$ is a double zero of the second-order elliptic function $\wp$. As its poles are also double, the Weierstrass function $\wp$ itself has meromorphic square-roots (by the Weierstrass product theorem, for instance). It also follows that if $g_3 \neq 0$ then none of the midpoint constants vanishes, so that $\wp$ has simple zeros and no meromorphic square-roots.

\medbreak 

Neville shifts $\wp$ by the midpoint constants and considers the three functions $\wp - e_p$ as $p$ runs over $\{ 1, 2, 3 \}$. By design, each of these second-order elliptic functions has double zeros on the corresponding midpoint lattice $\omega_p + \Lambda$ and so has two meromorphic square-roots; Neville (though with an ingenious change of notation, which we recommend) defines the primitive function $J_p$ to be the meromorphic square-root of $\wp - e_p$ that satisfies $z J_p (z) \rightarrow 1$ as $z \rightarrow 0$. Of course, our observation calls for no correction to any of this: it is simply the case that if $g_3 = 0$ then one of the midpoint constants is actually zero and need not be subtracted; the corresponding primitive function is then naturally preferred. 

\medbreak 

To take a particularly straightforward example, let $\omega_1 = 1$ and $\omega_2 = i$ so that $\omega_3 = -1 -i$ and 
$$\Lambda = \{ 2m + 2n i : m, n \in \mathbb{Z} \}$$
is the lattice of (even) Gaussian integers; the union $\frac{1}{2} \Lambda$ of $\Lambda$ and its three midpoint lattices is the full lattice of Gaussian integers. As multiplication by $i$ leaves $\Lambda$ invariant, 
$$\sum_{0 \neq \lambda \in \Lambda} \lambda^{-6} = \sum_{0 \neq \lambda \in \Lambda} (i \lambda)^{-6} = i^{-6} \sum_{0 \neq \lambda \in \Lambda} \lambda^{-6} = - \sum_{0 \neq \lambda \in \Lambda} \lambda^{-6}$$
whence
$$g_3(\Lambda) = 0$$ 
and a similar calculation reveals that 
$$\wp (i z) = - \wp (z).$$
It follows that $e_3 = \wp (\omega_3) = 0$: indeed, $\wp(\omega_2) = \wp(i) = - \wp(1) = -\wp(\omega_1)$ so that $e_1 + e_2 = 0$ while $e_1 + e_2 + e_3 = 0$ in any case; of course, a direct computation is also possible. For this lattice, the Weierstrass function $\wp$ has global meromorphic square-roots, namely $J_3$ and $-J_3$.  It may be checked that the identity $\wp(i z) = - \wp(z)$ implies that  $i J_3(i z) = J_3 (z)$; it may also be checked that the same symmetry interchanges the other primitive elliptic functions $J_1$ and $J_2$ in the sense $J_2 (z) =  i J_1 (i z)$ and $J_1 (z) = i J_2 (i z)$. 

\medbreak 

To summarize the general situation: if the invariant $g_3(\Lambda)$ vanishes, then one of the midpoint constants vanishes, naturally singling out the corresponding midpoint lattice along with the corresponding primitive elliptic function, which is a meromorphic square-root of the Weierstrass function $\wp_{\Lambda}$ itself; if $g_3 (\Lambda)$ does not vanish, then $\wp_{\Lambda}$ lacks meromorphic square-roots. To put a part of this another way, the invariant $g_3(\Lambda)$ is the obstruction to the existence of a global meromorphic square-root of $\wp_{\Lambda}$. 

\medbreak 

Incidentally, an obstruction-theoretic significance also attaches to the invariant 
$$g_2 = g_2 (\Lambda) =60 \sum_{0 \neq \lambda \in \Lambda} \lambda^{-4}$$
which satisfies 
$$e_2 e_3 + e_3 e_1 + e_1 e_2 = - g_2/4.$$
Let $\zeta_4$ be the fourth-order Eisenstein function defined by 
$$\zeta_4 (z) = \sum_{\lambda \in \Lambda} (z - \lambda)^{-4}$$
so that 
$$\zeta_4 = \frac{1}{6} \wp '' = \wp^2 - \frac{1}{12} g_2.$$
Evidently, if $g_2$ is zero then $\zeta_4$ has the functions $\pm \wp$ as meromorphic square-roots. Assume instead that $g_2$ is nonzero and write $c^2 = \frac{1}{12} g_2$: if $\zeta_4 = \wp^2 - c^2 = (\wp - c)(\wp + c)$ is a square then its zeros must be double, so that $c$ and $-c$ are midpoint constants; as the three midpoint constants have zero sum, they are $\pm c$ and $0$, whence 
$$- \frac{1}{4} g_2 = e_2 e_3 + e_3 e_1 + e_1 e_2 = - c^2 = - \frac{1}{12} g_2$$
and therefore $g_2$ is zero, contrary to assumption. In short, $\zeta_4$ admits global meromorphic square-roots precisely when the invariant $g_2$ vanishes. 

\medbreak

\bigbreak

\begin{center} 
{\small R}{\footnotesize EFERENCES}
\end{center} 
\medbreak 

[1] E.H. Neville, {\it Jacobian Elliptic Functions}, Oxford University Press (1944). 

\medbreak 

[2] E.H. Neville, {\it Elliptic Functions: A Primer}, Pergamon Press (1971).

\end{document}